\documentclass[draft]{amsart}
\usepackage{amssymb} 
\newtheorem{theorem}{Theorem}[section]

\newtheorem{s}[theorem]{S}
\newtheorem{dif}[theorem]{D}
\newtheorem{g}[theorem]{G}

\theoremstyle{definition}

\theoremstyle{approach}

\numberwithin{equation}{section}

\begin{document}
\setcounter{page}{1}
\title[Fr\'echet algebras in abstract harmonic analysis]{Fr\'echet algebras in abstract harmonic analysis}
\author[Z. Alimohammadi and A. Rejali]{Z. Alimohammadi and A. Rejali}



\begin{abstract}
We provide a survey of the similarities and differences between Banach and Fr\'echet algebras including some known results and examples.
We also collect some important generalizations in abstract harmonic analysis; for example, the generalization of the concepts of vector-valued Lipschitz algebras, abstract Segal algebras, Arens regularity, amenability, weak amenability, ideal amenability, etc.
\end{abstract}

\maketitle \setcounter{section}{-1}

\section{\bf Introduction}\label{s0}

The class of Fr\'echet algebras which is an important class of locally convex algebras has been widely studied by many authors. For a full understanding of Fr\'echet algebras, one may refer to \cite{Goldmann,Hel}. In the class of Banach algebras, there are many concepts which were generalized to the Fr\'echet case. Examples of these concepts are as follows.

Biprojective Banach algebras were studied by several authors, notably A. Ya. Helemskii \cite{Hel01} and Yu. V. Selivanov \cite{Se0}.
A number of papers appeared in the literature concerned with extending the notion of biprojectivity and its results to the Fr\'echet algebras; see for example \cite{Pir0,Se1}.

Flat cyclic Fr\'echet modules were investigated by A. Yu. Pirkovskii \cite{Pir}. 
Let $\mathcal{A}$ be a Fr\'echet algebra and $I$ be a closed left ideal in $\mathcal{A}$. Pirkovskii showed that there exists a necessary and sufficient condition for a cyclic Fr\'echet $\mathcal{A}$-module $\mathcal{A}^{\sharp}/I$ to be strictly flat. He also generalized a number of characterizations of amenable Banach algebras to the Fr\'echet algebras (ibid.). In 2008, P. Lawson and C. J. Read \cite{LR} introduced and studied the notions of approximate amenability and approximate contractibility of Fr\'echet algebras. 
Recently, the concepts of $\varphi$-amenability and character amenability, $\varphi$-contractibility and character contractibility, weak amenability and Ideal amenability for Fr\'echet algebras have investigated in several papers; see \cite{ARR1,AR,ARR4,RR}. 

T. G. Honary \cite{Honary1} gave an extension of Johnson's uniqueness of norm theorem. He showed that if $\mathcal{A}$ and $\mathcal{B}$ are Fr\'echet algebras, $\mathcal{B}$ is semisimple and $T:\mathcal{A}\rightarrow\mathcal{B}$ is a surjective homomorphism with a certain condition, then $T$ is continuous. In special case, if $\mathcal{A}$ is a Banach algebra, then every epimorphism $T:\mathcal{A}\rightarrow\mathcal{B}$ is automatically continuous. Hence, every semisimple Banach algebra has a unique topology as a Fr\'echet algebra. Honary then presented a partial answer to the well-known Michael's problem posed in \cite{Michael}.

Let $(X,d)$ be a metric space and $E$ be a Banach algebra over $\mathbb{C}$. Take $\alpha\in\mathbb{R}$ with $\alpha>0$. E. Biyabani and A. Rejali \cite{BR} studied the structure and properties of vector-valued Lipschitz algebras $Lip_{\alpha}(X,E)$ and $lip_{\alpha}(X,E)$ of order $\alpha$. Later, these concepts were generalized to the class of Fr\'echet algebras by Rejali 
joint with Ranjbari \cite{RR0}. They also studied the concept of ideal amenability for vector-valued Fr\'echet Lipschitz algebras (ibid.). Furthermore, Rejali et al. \cite{ARR2,ARR3} introduced and studied the notions of Segal and semisimple Segal Fr\'echet algebras.
They showed that every continuous linear left multiplier of a Fr\'echet algebra $(\mathcal{B},q_n)_{n\in\mathbb{N}}$ is also a continuous linear left multiplier of any Segal Fr\'echet algebra $(\mathcal{A},p_{\ell})_{\ell\in\mathbb{N}}$ in $(\mathcal{B},q_n)_{n\in\mathbb{N}}$. Moreover, they proved that if $\mathcal{B}$ is a commutative Fr\'echet $Q$-algebra with an approximate identity, then the spaces of all modular maximal closed ideals of $\mathcal{B}$ and any Segal Fr\'echet algebra $\mathcal{A}$ in $\mathcal{B}$ are homeomorphic. Particularly, $\mathcal{B}$ is semisimple if and only if $\mathcal{A}$ is semisimple \cite{ARR3}.
In \cite{AA}, we joint with Abtahi introduced and studied the strict, uniform, and compact-open locally convex topologies on an arbitrary algebra $\mathcal{B}$, by the fundamental system of seminorms of a locally convex subalgebra $(\mathcal{A},p_{\alpha})$. In particular, we studied the strict topology on the multiplier algebra $M(\mathcal{A})$ of a Fr\'echet algebra $\mathcal{A}$, and we showed that $\mathcal{A}$ is a $M(\mathcal{A})$-Segal algebra. 

In the present paper we give a list of the significant results which have been obtained by various authors. Some of the results in this paper were included in the first author's Ph.D. thesis, which was supported by the center of excellence for mathematics at the University of Isfahan.

\section{\bf Preliminaries}

In this section we give some basic definitions and frameworks related to locally convex spaces
which will be required throughout the paper; see \mbox{\cite{Goldmann,Hel,Meise}} for more information.
 In this paper, all vector spaces and algebras are assumed to be over the field $\mathbb{C}$ of complex numbers.

 A locally convex space $E$ is a topological vector space
in which the origin has a local base of absolutely convex
absorbent sets. Throughout the paper, all locally convex spaces
are assumed to be Hausdorff. Every locally convex Hausdorff space
$E$ has a fundamental system of seminorms $(p_{\alpha})_{\alpha\in
\Lambda}$, or equivalently a family of continuous seminorms
satisfying the following properties:
\begin{enumerate}
\item[(i)] For every $x\in E$ with $x\neq 0$, there exists an
${\alpha}\in \Lambda$ with $p_{\alpha}(x)>0$; \item[(ii)] For all
${\alpha,\beta}\in \Lambda$, there exist ${\gamma}\in \Lambda$ and $C>0$ such
that $$\max(p_{\alpha}(x),p_{\beta}(x))\leq C
p_{\gamma}(x)\;\;\;\;\;\;\;\;(x\in E).$$
\end{enumerate}
We denote by $(E,p_{\alpha})_{\alpha\in
\Lambda}$, the locally convex space with the fundamental system of seminorms $(p_{\alpha})_{\alpha\in
\Lambda}$.

Suppose that $(E,p_{\alpha})_{\alpha\in
\Lambda}$ and $(F,q_{\beta})_{\beta \in B}$ are locally convex spaces. For every linear
mapping $T: E\rightarrow F$, the following assertions are
equivalent:
\begin{enumerate}
\item[(i)]  $T\in \mathcal{L}(E,F)$, i.e. $T$ is continuous;
\item[(ii)] $T$ is continuous at $0$; 
\item[(iii)] For each $\beta\in B$ there
exist an $\alpha\in A$ and $C>0$, such that $$q_{\beta}(T(x))\leq
Cp_{\alpha}(x)\;\;\;\;\;\;\;\;\;\;\;(x\in E).$$
\end{enumerate}

Following \cite{Meise}, a complete metrizable locally convex space is called a Fr\'echet
space. Indeed, a Fr\'echet space is a locally convex space which
has a countable fundamental system of seminorms
$(p_{\ell})_{\ell\in\Bbb N}$. It is worth noting to remark that due to Conway \cite{Conway}, Fr\'echet spaces are not assumed to be locally convex. These two different definitions make some basic differences. Consider $L^p(0,1)$ ($0<p<1$), the space of all equivalence classes of measurable functions \mbox{$f:(0,1)\rightarrow\mathbb{R}$} such that $\int_0^1|f(x)|^p\mathrm{d}x<\infty$. Under the translation metric 
$$
d(f,g)=\int_0^1|(f-g)x|^p\mathrm{d}x<\infty\;\;\;\;\;\;\;\;\;\;\;(f,g\in L^p(0,1)),
$$
$L^p(0,1)$ is a complete metric space. Therefore, according to the definition presented in \cite{Conway}, $L^p(0,1)$ is a Fr\'echet space. However, $L^p(0,1)$ has only one nonempty open convex set, namely itself. It follows that due to the definition given in \cite{Meise}, $L^p(0,1)$ is not a Fr\'echet space. In this paper, our definition of a Fr\'echet space is based on that given in \cite{Meise}. 

A Fr\'echet space $\mathcal{A}$ is called quasinormable if for each zero neighborhood $U$ there is a zero neighborhood $V$ so that for every $\varepsilon>0$ there exists a bounded set $B$ in $\mathcal{A}$ such that $V\subseteq B+\varepsilon U$.

It should be noted that, by \cite[Proposition 25.7]{Meise}, the strong dual space of every Fr\'echet space $\mathcal{A}$ is a (DF)-space which is a type of locally convex topological vector spaces; see \cite[page 297]{Meise} for details. We denote the strong dual of $\mathcal{A}$ by $\mathcal{A}^*$.

Following \cite[Corollary 24.29]{Meise}, every Fr\'echet space $\mathcal{A}$ has a web $\{C_{n_1,\cdots,n_k}\}$, that is a family $C_{n_1,\cdots,n_k}$, $n_1,\cdots,n_k\in\mathbb{N}$, $k\in\mathbb{N}$, of absolutely convex subsets of $\mathcal{A}$ with the following properties:
\begin{enumerate}
\item[(i)]
$\cup_{n=1}^{\infty}C_n=\mathcal{A}$;
\item[(ii)]
$\cup_{n=1}^{\infty}C_{n_1,\cdots,n_k,n}=C_{n_1,\cdots,n_k}$ for all $n_1,\cdots,n_k\in\mathbb{N}$ and all $k\in\mathbb{N}$;
\item[(iii)]
for each sequence $(n_k)_{k\in\mathbb{N}}$ in $\mathbb{N}$ there exists a sequence $(\lambda_k)_{k\in\mathbb{N}}$ in $(0,\infty)$, so that for every sequence $(x_k)_{k\in\mathbb{N}}$ in $\mathcal{A}$ with $x_k\in C_{n_1,\cdots,n_k,n}$ for all $k\in\mathbb{N}$ the series $\sum_{k=1}^{\infty}\lambda_kx_k$ converges in $\mathcal{A}$.
\end{enumerate}

A barrel set in a locally convex space $E$ is a subset which is convex, balanced, absorbing and closed. Moreover, $E$ is called barrelled if each barrel in $E$ is a zero neighborhood. 
By \cite[Proposition 23.23]{Meise}, every Fr\'echet space is barrelled.

Let $\Lambda$ be an index set directed under a (reflexive, transitive, anti-symmetric) relation $\prec$, let  $(E_{\alpha})_{\alpha\in\Lambda}$ be a family of  locally
convex spaces, and denote, for $\alpha\prec\beta$, by $\tau^{\beta}_{\alpha}$ a continuous linear map of $E_{\beta}$ into $E_{\alpha}$. Consider
$$
E=\big\{ (x_{\alpha})_{\alpha\in\Lambda}\in\prod_{\alpha\in\Lambda} E_{\alpha}:
\tau^{\beta}_{\alpha}(x_{\beta})=x_{\alpha}\;\text{whenever}\;\alpha\prec\beta\big\}.
$$
$E$ is called the projective limit of the family $(E_{\alpha})_{\alpha\in\Lambda}$ with respect to the mappings $\tau^{\beta}_{\alpha}$ ($\alpha,\beta\in\Lambda,\;\alpha\prec\beta$), and denoted by $\underleftarrow{\lim}(E_{\alpha},\tau^{\beta}_{\alpha})$. 
The topology of $E$ is the projective topology with respect to the family $\{(E_{\alpha},t_{\alpha},\tau_{\alpha})\}_{\alpha\in\Lambda}$, where for each $\alpha\in\Lambda$, $t_{\alpha}$ is the topology of $E_{\alpha}$ and $\tau_{\alpha}:E\rightarrow E_{\alpha}$ is the canonical map. 
The projective limit $E$ is called reduced if for each $\alpha\in\Lambda$, $\tau_{\alpha}$ has dense range. By applying \cite[II.5.4]{Sch}, every Fr\'echet space is isomorphic to a projective limit of a sequence of Banach spaces.

By a similar way, one can introduce the inductive limit of the family $(E_{\alpha})_{\alpha\in\Lambda}$. See also \cite[II.6.3]{Sch} and \cite[section 24]{Meise} for more information about the notion of inductive topology.

A locally convex space $E$ is called bornological (resp. ultra-bornological)
if its topology is the inductive topology corresponding to a system $(E_{\alpha}\rightarrow E)_{\alpha\in\Lambda}$
of normed (resp. Banach) spaces. Clearly, every ultra-bornological space is bornological. Also, every Fr\'echet space is ultra-bornological \cite[Remark 24.13]{Meise}.

A topological algebra ${\mathcal A}$ is an algebra, which is a
topological vector space and the multiplication ${\mathcal
A}\times{\mathcal A}\rightarrow {\mathcal A}$ $((a,b)\mapsto
ab)$ is separately continuous.
Moreover, ${\mathcal A}$ is called a locally convex algebra if ${\mathcal A}$ in addition is a locally
convex space. 

Let $(\mathcal{A},p_{\alpha})_{\alpha\in\Lambda}$ be a locally convex algebra. A left approximate identity (left a.i. for short) in $\mathcal{A}$ is a net $(e_i)_{i\in I}$ in $\mathcal{A}$, such that 
$$
p_{\alpha}(e_ia-a)\longrightarrow_i0\;\;\;\;\;\;\;\;\;\;\;(a\in\mathcal{A})
$$
for each $\alpha\in\Lambda$. The left a.i. is bounded (left b.a.i. for short) if the set $\{e_i:i\in I\}$ is bounded. It is called uniformly bounded (left u.b.a.i. for short) if 
$$
\sup\{p_{\alpha}(e_i):i\in I,\;\alpha\in\Lambda\}<\infty.
$$
``Right a.i.", ``right b.a.i.", and ``right u.b.a.i." are defined similarly. Furthermore, for two-sided approximate identities, we use of the symbols ``a.i.", ``b.a.i.", and ``u.b.a.i.".

Following \cite{Goldmann}, a complete topological algebra is a Fr\'echet algebra if its topology is produced by a countable family of increasing submultiplicative seminorms.

Let $\mathcal{A}$ be a Fr\'echet algebra. A left Fr\'echet $\mathcal{A}$-module is a Fr\'echet space $E$ together with the structure of an $\mathcal{A}$-module such that the bilinear mapping
$$
\mathcal{A}\times E\rightarrow E,\;\;\;\;(a,x)\mapsto ax
$$
is continuous.
Similarly, right Fr\'echet $\mathcal{A}$-modules and Fr\'echet $\mathcal{A}$-bimodules are defined.
Suppose that $E$ and $F$ are left Fr\'echet $\mathcal{A}$-modules (resp. right Fr\'echet $\mathcal{A}$-modules, Fr\'echet $\mathcal{A}$-bimodules). In this paper, the space of all continuous $\mathcal{A}$-module morphisms from $E$ to $F$ is denoted by $_{\mathcal{A}}\mathcal{L}(E,F)$ (resp. $\mathcal{L}_{\mathcal{A}}(E,F)$, $_{\mathcal{A}}\mathcal{L}_{\mathcal{A}}(E,F)$).

\section{\bf Comparison between Banach and Fr\'echet algebras}

In this section, we compare two Banach and Fr\'echet algebras and provide a list of important similarities and differences between them. We use the symbol ``$\Box$" after each result to separate them from one another in the list.

\subsection{Similarities}\label{2}

Following results are marked ``S".

\begin{s}
Due to \mbox{\cite[Proposition 22.12]{Meise}}, there exists a generalization of the Hahn-Banach theorem and its consequences for locally convex spaces.  Let $\mathcal{A}$ be a locally convex space, $B$ be a subspace of $\mathcal{A}$, and $p$ be a continuous seminorm on $\mathcal{A}$. Then, the following statements hold.
\begin{enumerate}
\item[(i)]
For each $f\in B^*$ there exists an $F\in\mathcal{A}^*$ such that $F|_{B}=f$.
\item[(ii)]
For each $a\in\mathcal{A}$ there exists an $f\in\mathcal{A}^*$ such that $f(a)=p(a)$ and $|f|\leq p$.
\item[(iii)]
For each $a\in\mathcal{A}$ with $a\not=0$ there exists an $f\in\mathcal{A}^*$ such that $f(a)\not=0$. 
\end{enumerate}
{\bf Example.}
Consider the space $\mathbb{C}^{\mathbb{N}}$ with the fundamental system of seminorms
$$
p_n(x)=\max_{1\leq m\leq n}|x_m|\;\;\;\;\;\;\;\;\;\;\;(x=(x_n)_{n\in\mathbb{N}}\in\mathbb{C}^{\mathbb{N}}).
$$
Let the subspace $B$ consist of all convergent complex-valued sequences. By arguments above, the map $T:B\rightarrow\mathbb{C}$, defined by
$$
T(x)=\lim_{n}x_n\;\;\;\;\;\;\;\;\;\;\;(x=(x_n)_{n\in\mathbb{N}}\in B),
$$
has a continuous linear extension to $\mathbb{C}^{\mathbb{N}}$. This extension is denoted by LIM; see \cite[Theorem 1]{Hajdukovic}, \cite[III.9]{Sherbert} and also \cite[section 4]{RR0} for details. $\Box$
\end{s}

Before the next result, let us recall some notation.
Let $\mathcal{A}$ be a locally convex space. The polars, of nonempty sets $M\subseteq\mathcal{A}$ and $N\subseteq\mathcal{A}^*$, are defined by
\begin{eqnarray*}
M^{\circ}&=&\big{\{}f\in\mathcal{A}^*: |f(a)|\leq 1\; \text{for all}\; a\in M\big{\}}\subseteq \mathcal{A}^*,\;\text{and} \\
^{\circ}N&=&\big{\{}a\in\mathcal{A}: |f(a)|\leq 1\; \text{for all}\; f\in N\big{\}}\subseteq \mathcal{A}.
\end{eqnarray*}

\begin{s}
Following \mbox{\cite[Theorem 22.13]{Meise}}, the bipolar theorem is valid in the class of locally convex spaces. 
Let $\mathcal{A}$ be a locally convex space and $B$ be an absolutely convex subset of $\mathcal{A}$. Then, $^{\circ}(B^{\circ})=\overline{B}$. $\Box$
\end{s}

\begin{s}
We know that if $\mathcal{A}$ is a normed space, then by applying Banach-Alaoglu theorem, the closed unit ball of $\mathcal{A}^*$ is weak$^*$-compact. The Alaoglu-Bourbaki theorem is a generalization by Bourbaki to dual topologies on locally convex spaces. Let $\mathcal{A}$ be a locally convex space. Then, the polar $U^{\circ}$ of any neighborhood $U$ in $\mathcal{A}$ is weak$^*$-compact; see \mbox{\cite[Theorem 23.5]{Meise}}. $\Box$
\end{s}

\begin{s}\label{s4}
Let $\mathcal{A}$ and $\mathcal{B}$ be Banach spaces. Due to the open mapping theorem,
if the map $T:\mathcal{A}\rightarrow\mathcal{B}$ is continuous, linear and surjective, then $T$ is open.
A similar result holds for Fr\'echet spaces instead of Banach spaces.
Following \mbox{\cite[Theorem 24.30]{Meise}}, the open mapping theorem was generalized to the locally convex spaces. Let $\mathcal{A}$ and $\mathcal{B}$ be locally convex spaces. If $\mathcal{A}$ has a web and $\mathcal{B}$ is ultra-bornological, then every continuous, linear, surjective map $T:\mathcal{A}\rightarrow\mathcal{B}$ is open\mbox{. $\Box$}
\end{s}

\begin{s}
Consider the continuous map $f:\mathcal{A}\rightarrow\mathcal{B}$ between two metric spaces $\mathcal{A}$ and $\mathcal{B}$. Following \cite{Meise}, 
$$ 
\mathcal{G}(f)=\big{\{}(a,f(a)):a\in\mathcal{A}\big{\}}\subseteq\mathcal{A}\times\mathcal{B},
$$
denote the graph of $f$, is closed. If $\mathcal{A}$ and $\mathcal{B}$ are also complete, then the converse is true. 
 In general, the closed graph theorem holds in the class of locally convex spaces. Suppose that $\mathcal{A}$ and $\mathcal{B}$ are locally convex spaces. Exactly like \mbox{S \ref{s4}}, if $\mathcal{A}$ has a web and $\mathcal{B}$ is ultra-bornological, then every linear map $T:\mathcal{A}\rightarrow\mathcal{B}$ with a closed graph is continuous \mbox{\cite[Theorem 24.31]{Meise}}.
In particular, if $\mathcal{A}$ and $\mathcal{B}$ are Fr\'echet spaces, then a linear mapping of $\mathcal{A}$ into $\mathcal{B}$ is continuous if and only if its graph is closed in $\mathcal{A}\times\mathcal{B}$. $\Box$ 
\end{s}

\begin{s}
In addition to the open mapping theorem and closed graph theorem, there is an extension of the Banach-Steinhaus theorem, known as uniform boundedness principle. Let $\mathcal{A}$ and $\mathcal{B}$ be locally convex spaces and $\mathcal{A}$ be barrelled.
Further, let $H\subseteq \mathcal{L}(\mathcal{A},\mathcal{B})$ be pointwise bounded, i.e., for each $a\in\mathcal{A}$, $\{Ta:T\in H\}$ is bounded in $\mathcal{B}$. Then, $H$
is equicontinuous. In other words, for each zero neighborhood $V$ in $\mathcal{B}$ there exists a zero neighborhood $U$ in $\mathcal{A}$ with $T(U)\subseteq V$ for all $T\in H$; see \mbox{\cite[Theorem III.\S 4.2.1]{Bourbaki}} and \mbox{\cite[Proposition 23.27]{Meise}}. \\
{\bf Corollary.}
Let $\mathcal{A}$ and $\mathcal{B}$ be locally convex spaces, $\mathcal{A}$ be barrelled and $(T_n)_{n\in\mathbb{N}}$ be a sequence in $\mathcal{L}(\mathcal{A},\mathcal{B})$. Also, let $(T_na)_{n\in\mathbb{N}}$ be convergent in $\mathcal{B}$ for each $a\in\mathcal{A}$. Then, the mapping $T:\mathcal{A}\rightarrow\mathcal{B}$, defined by
$$
Ta=\lim_{n\rightarrow\infty}T_na\;\;\;\;\;\;\;\;\;\;\;(a\in\mathcal{A}),
$$ 
is continuous \mbox{\cite[page 274]{Meise}}. $\Box$
\end{s}

\begin{s}
Let $\mathcal{A}$ be a Fr\'echet algebra, $E$ a right Fr\'echet $\mathcal{A}$-module, and $F$ a left Fr\'echet $\mathcal{A}$-module. By applying \mbox{\cite[Theorem 2]{Smith}}, the completed projective tensor product of $E$ and $F$, denoted by $E\hat{\otimes}F$, is again a Fr\'echet algebra. Now, denote by $E\hat{\otimes}_{\mathcal{A}}F$ the completion of the quotient $(E\hat{\otimes}F)/N$, where  
$N\subseteq E\hat{\otimes}F$ is the closed linear span of all elements of the form 
$$
x\cdot a\otimes y-x\otimes a\cdot y\;\;\;\;\;\;\;\;\;\;\;(x\in E,\;y\in F,\;a\in\mathcal{A}).
$$ 
Following \cite[Proposition 3.2]{Pir}, there exists a vector space isomorphism from $(E\hat{\otimes}_{\mathcal{A}}F)^*$ to $\mathcal{L}_{\mathcal{A}}(E,F^*)$. $\Box$ 
\end{s}

Now, let $\mathcal{A}$ be a commutative Fr\'echet algebra. For each $a\in\mathcal{A}$, let $\hat{a}$ be the Gelfand transform of $a$ defined by $\hat{a}(\varphi):=\varphi(a)$, where $\varphi$ is an arbitrary nonzero complex-valued algebra homomorphism on $\mathcal{A}$.
Note that every complex-valued algebra homomorphism on a Banach algebra is automatically continuous. However, the question of whether each homomorphism on a Fr\'echet algebra is necessarily continuous was posed by E. A. Michael \cite{Michael}. This question has not been solved. But, only partial answers have been obtained so far; for example see \cite{Honary1} as mentioned in the introduction. See also \mbox{\cite[chapter 10]{Goldmann}} for details.
In this paper, the set of all nonzero complex-valued and continuous algebra homomorphisms on $\mathcal{A}$ is denoted by $\Delta(\mathcal{A})$. Also, the topology on $\Delta(\mathcal{A})$ is the Gelfand topology,  that is the coarsest topology such that all 
Gelfand transforms are continuous functions on $\Delta(\mathcal{A})$.  

\begin{s}\label{spectrum}
Let $\mathcal{A}$ be a commutative unital Fr\'echet algebra.
Following \mbox{\cite[3.2]{Goldmann}}, $\Delta(\mathcal{A})$ is never empty. Therefore, the Gelfand-Mazur theorem holds for Fr\'echet algebras. In other words, if every nonzero element in $\mathcal{A}$ is invertible, then $\mathcal{A}=\mathbb{C}$. $\Box$ 
\end{s}

Let $\mathcal{B}$ be a Banach algebra with a left b.a.i. (not necessarily u.b.a.i.) and let $E$ be a left Banach $\mathcal{B}$-module. Due to the Cohen factorization theorem, 
for given $\varepsilon>0$ and $x\in \overline{\mathcal{B}\cdot E}$, there are $a\in\mathcal{B}$ and $y\in E$ such that $x=a\cdot y$ and $\|x-y\|<\varepsilon$. Moreover, $\|a\|<M$, for some $M\in\mathbb{R}^{+}$; see \cite[page 14]{Dales}.
We mention next an important generalization by Voigt of the Cohen factorization theorem; see \cite{Voigt} for details.

\begin{s}
Let $(\mathcal{A},p_{\ell})$ be a Fr\'echet algebra and let $(e_i)_{i\in I}$ be a left u.b.a.i. for $\mathcal{A}$. 
Suppose that $E$ is a left Fr\'echet $\mathcal{A}$-module and $\mathtt{B}\subseteq E$ is bounded.
Let us recall $aco(\mathtt{B})$ denote the absolutely convex hull of $\mathtt{B}$, that is 
$$
aco(\mathtt{B})=\big{\{}\sum_{i=1}^{n}\lambda_ix_i:x_i\in \mathtt{B},\;\lambda_i\in\mathbb{C},\;\sum_{i=1}^{n}|\lambda_i|\leq1\big{\}}.
$$
Set $\mathfrak{B}:=\overline{aco}(\mathtt{B})$ and $E_{\mathfrak{B}}:=(lin(\mathfrak{B}),p_{\mathfrak{B}})$, where $p_{\mathfrak{B}}$ is the gauge $($or Minkowski functional$)$ of $\mathfrak{B}$, i.e., the non-negative real function
$$
x\mapsto p_{\mathfrak{B}}(x)=\inf\{t>0:x\in t\mathfrak{B}\}\;\;\;\;\;\;\;\;\;\;\;(x\in\mathtt{B}).
$$
By \cite[Theorem 3]{Voigt}, if $e_ix\longrightarrow_i x$ uniformly for $x\in \mathtt{B}$, then $(\mathcal{A},E)$ has the \mbox{$\{\mathtt{B}\}$-strong} left factorization property (hereafter abbreviated $\{\mathtt{B}\}$-SLFP). In other words,
for every zero neighborhood $U$ in $E$, there exist $a\in\mathcal{A}$ and a continuous linear mapping $T:E_{\mathfrak{B}}\rightarrow E$ such that
\begin{enumerate}
\item[(i)]
$a(Tx)=x$, for all $x\in \mathtt{B}$;
\item[(ii)]
$Tx$ belongs to the closed $\mathcal{A}$-submodule of $E$ generated by $x$, for all $x\in \mathtt{B}$;
\item[(iii)]
$Tx-x\in U$, for all $x\in \mathtt{B}$;
\item[(iv)]
$T|_{\mathtt{B}}:\mathtt{B}\rightarrow E$ is continuous.
\end{enumerate}
Furthermore, $a\in\mathcal{A}$ can be chosen to satisfy
\begin{enumerate}
\item[(v)]
if $\ell\in\mathbb{N}$ is such that $\sup\{p_{\ell}(e_i): i\in I\}\leq 1$, then $p_{\ell}(a)\leq 1$.
\end{enumerate}
{\bf Remark.}
Let $\mathcal{B}$ be a collection of bounded (resp. compact) subsets of $E$. Following \cite{Voigt}, we say that  $(\mathcal{A},E)$ has the $\mathcal{B}$-SLFP if for each $\mathtt{B}\in\mathcal{B}$, $(\mathcal{A},E)$ has the $\{\mathtt{B}\}$-SLFP. In particular, if $\mathcal{B}$ is the collection of bounded (resp. compact) subsets of $\mathcal{A}$, then this property is called the bounded (resp. compact) SLFP. 

Note that by \cite[Corollary 5]{Voigt}, if $\mathcal{A}$ has a left u.b.a.i., then $\mathcal{A}$ has the compact SLFP. However, factorization properties are proved for certain Fr\'echet algebras of differentiable functions which do not possess a u.b.a.i.; see \cite{Voigt2}. $\Box$
\end{s}

\subsection{Differences}\label{3}

Following results are marked ``D".

\begin{dif}
Every Banach algebra (space) is a Fr\'echet algebra (space), but in general the converse is not true. 
Let $(\mathcal{A}_n,\|\cdot\|_n)_{n\in\mathbb{N}}$ be a sequence of normed spaces and  $\mathcal{A}=\prod_{n\in\mathbb{N}}\mathcal{A}_n$. Define the metric $d$ on $\mathcal{A}$ by
$$
d(a,b):=\sum_{n=1}^{\infty}\frac{\|a_n-b_n\|_n}{2^n\big{(}1+\|a_n-b_n\|_n\big{)}}\;\;\;\;\;\;\;\big{(}a=(a_n)_{n\in\mathbb{N}}, b=(b_n)_{n\in\mathbb{N}}\in\mathcal{A}\big{)}.
$$
Following \cite[Lemma 5.17]{Meise}, $(\mathcal{A},d)$ is a locally convex metric linear space. Also,  if each $(\mathcal{A}_n,\|\cdot\|_n)$ is complete, then $(\mathcal{A},d)$ is a Fr\'echet space. However, $(\mathcal{A},d)$ is not a Banach space if $\mathcal{A}_n\not=\{0\}$ for infinitely many $n\in\mathbb{N}$. For instance, the space $\mathbb{C}^{\mathbb{N}}$ of all sequences in $\mathbb{C}$ with the metric 
$$
d(x,y)=\sum_{n=1}^{\infty}\frac{|x_n-y_n|_n}{2^n\big{(}1+|x_n-y_n|_n\big{)}}\;\;\;\;\;\;\;\big{(}x=(x_n)_{n\in\mathbb{N}}, y=(y_n)_{n\in\mathbb{N}}\in\mathbb{C}^{\mathbb{N}}\big{)},
$$
is a Fr\'echet space which is not a Banach space; see \mbox{\cite[Examples 5.18]{Meise}}. $\Box$
\end{dif}

\begin{dif}\label{3.222}
The dual space of a Banach space is also a Banach space. Moreover, if $\mathcal{A}$ is a Banach algebra, then the module actions of $\mathcal{A}$ are jointly continuous on $\mathcal{A}^*$. However, this result is not valid for a Fr\'echet case.
Let $C^{\infty}(\mathbb{R})$, denote the space of all infinitely differentiable functions on $\mathbb{R}$. Following \cite{Tay}, $C^{\infty}(\mathbb{R})$ is a Fr\'echet algebra with the pointwise product and the family of countable seminorms $\|\cdot\|_n$, defined by
$$
\|f\|_n:=\sum_{k=0}^{n}\frac{1}{k!}\|f^{(k)}\|_{\infty,n}\;\;\;\;\;\;\;\;(f\in C^{\infty}(\mathbb{R})).
$$ 
The module action of $C^{\infty}(\mathbb{R})$ is not  jointly continuous on $C^{\infty}(\mathbb{R})^*$. Therefore, $C^{\infty}(\mathbb{R})^*$ is not a Fr\'echet space. By arguments above, we immediately get the following. \\
{\bf Lemma.}
Let $\mathcal{A}$ be a metrizable locally convex space. The following statements are equivalent:
\begin{enumerate}
\item[(i)]
$\mathcal{A}$ is a Banach space;
\item[(ii)]
$\mathcal{A}$ is a Fr\'echet and (DF)-space.
\end{enumerate} 
{\bf Corollary.}
Let $\mathcal{A}$ be a Fr\'echet space. Then, $\mathcal{A}^*$ is a Fr\'echet space if and only if $\mathcal{A}$ is a Banach space. \\
See also \cite{Goldmann} and \cite{Meise}. $\Box$
\end{dif}

\begin{dif}
For Banach spaces $\mathcal{A}$ and $\mathcal{B}$, $\mathcal{L}(\mathcal{A},\mathcal{B})$ is also a Banach space; see \mbox{\cite[Proposition 5.6]{Meise}}.  
However, by using \mbox{D \ref{3.222}} and \mbox{\cite[page 253]{Meise}}, $\mathcal{L}(\mathcal{A},\mathcal{B})$ is not necessarily a Fr\'echet space when $\mathcal{A}$ and $\mathcal{B}$ are Fr\'echet spaces. $\Box$
\end{dif}

Before proceeding to the next result, recall that a matrix $A=(a_{j,k})_{j,k\in\mathbb{N}}$ of non-negative numbers is called a K\"othe matrix if it satisfies the following conditions.
\begin{enumerate} 
\item[(i)]
For each $j\in\mathbb{N}$ there exists a $k\in\mathbb{N}$ with $a_{j,k}>0$.
\item[(ii)]
$a_{j,k}\leq a_{j,k+1}$ for all $j,k\in\mathbb{N}$.
\end{enumerate}
Now, for $1\leq p<\infty$ define
$$
\lambda^p(A):=\big{\{}x\in\mathbb{C}^{\mathbb{N}}: \|x\|_{k}:=\big{(}\sum_{j=1}^{\infty}|x_ja_{j,k}|^p\big{)}^{1/p}<\infty\;\text{for all $k\in\mathbb{N}$}\big{\}}.
$$
Also, for $p=\infty$ and $p=0$ define
\begin{eqnarray*}
\lambda^{\infty}(A)&:=&\big{\{}x\in\mathbb{C}^{\mathbb{N}}: \|x\|_{k}:=\sup_{j\in\mathbb{N}}|x_j|a_{j,k}<\infty\;\text{for all $k\in\mathbb{N}$}\big{\}},\;\text{and}\\
c_0(A)&:=&\big{\{}x\in\lambda^{\infty}(A): \lim_{j\rightarrow\infty}x_ja_{j,k}=0\;\text{for all $k\in\mathbb{N}$}\big{\}},
\end{eqnarray*}
respectively. Using the definition above, we can present in \mbox{D \ref{dis3.4d}} an example of \mbox{non-distinguished} Fr\'echet space due to Meise and Vogt \cite{Meise}.
Note that a metrizable locally convex space $\mathcal{A}$ is distinguished if its strong dual is a barrelled or bornological locally convex space. Consider the Fr\'echet space $\mathcal{A}$ and let $\mathcal{A}=\underleftarrow{\lim}(\mathcal{A}_n,\tau^m_n)_{n\in\mathbb{N}}$. By applying \cite[Remark 25.13]{Meise}, if $\mathcal{A}$ is distinguished, then $\mathcal{A}^*\cong\underrightarrow{\lim}(\mathcal{A}^*_n,(\tau^m_n)^*)_{n\in\mathbb{N}}$ as two locally convex spaces. In other words, the inductive topology coincides with the strong topology on $\mathcal{A}^*$ if and only if $\mathcal{A}$ is distinguished. This is an important advantage of distinguished spaces.  

\begin{dif}\label{dis3.4d}
Every Banach space is distinguished \cite[Corollary II.7.1]{Sch}. In general, however, a Fr\'echet space fails to be distinguished. For example, consider the K\"othe matrix $A=(a_{i,j;k})_{(i,j)\in\mathbb{N}^2,k\in\mathbb{N}}$ such that
\begin{enumerate}
\item[(i)] 
for each $i,j,k\in\mathbb{N}$ with $k\leq i$, $a_{i,j;k}=a_{i,j;1}>0$; 
\item[(ii)]
for each $m\in\mathbb{N}$, $\lim_{j\rightarrow\infty}a_{m,j;m}a^{-1}_{m,j;m+1}=0$. 
\end{enumerate}
By using \cite[Corollary 27.18]{Meise}, $\lambda^1(A)$ is a Fr\'echet space which is not distinguished. This example is a generalization of a classical counter-example of Grothendieck and K\"othe; see \cite[\S 31.7]{K1}. $\Box$
\end{dif}

\begin{dif}\label{d5}
Following \cite[Proposition 7.5]{Meise}, if $\mathcal{A}$ is a reflexive Banach space and $B$ is a closed subspace of $\mathcal{A}$, then $\mathcal{A}/B$ is reflexive. Also, every closed subspace of a reflexive Fr\'echet space is reflexive \cite[Proposition 23.26]{Meise}.
But, there are reflexive Fr\'echet spaces with non-reflexive quotients. Consider the K\"othe matrix $A=(a_{i,j;k})_{(i,j)\in\mathbb{N}^2,k\in\mathbb{N}}$, where
\begin{eqnarray*}
a_{i,j;k}:=
\begin{cases}
(ki)^k & \text{for $j<k,\;i\in\mathbb{N}$} \\
k^j & \text{for $j\geq k,\;i\in\mathbb{N}$}.
\end{cases}
\end{eqnarray*}
By \cite[Example 27.21]{Meise}, for each infinite subset $I$ of $\mathbb{N}^2$ and each $n\in\mathbb{N}$ there exists $k\in\mathbb{N}$ such that $\inf_{(i,j)\in I}a_{i,j;n}a^{-1}_{i,j;k}=0$. Thus, due to Dieudonn\'e-Gomes theorem \cite[Theorem 27.9]{Meise}, $\lambda^1(A)$ is reflexive. Since $a_{i,j;1}=1$ and \mbox{$a_{i,k;k}=a_{1,k;k}$}, $\lambda^1(A)$ has a quotient which is isomorphic to $\ell^1$; see \cite[Proposition 27.22]{Meise}. However, by applying \cite[Corollary 7.10]{Meise}, $\ell^1$ is non-reflexive\mbox{. $\Box$}
\end{dif}

For the next result, we need to remind the reader that a Montel space $\mathcal{A}$ is a quasi-barrelled space in which each bounded set is relatively compact. Let $A$ be the K\"othe matrix defined as in \mbox{D \ref{d5}}. By \cite[Theorem 27.9]{Meise}, $\lambda^p(A)$ is a Montel space for each $p\in [1,\infty]$. For this reason, $\lambda^p(A)^*$ is also a Montel space and so is reflexive; see \cite[Proposition 24.25]{Meise}. Following the proof of \mbox{\cite[Proposition 27.22]{Meise}}, the continuous linear map $Q:\lambda^p(A)\rightarrow\ell^p$ for $p\in [1,\infty)$ defined by
$$
Qx=(\sum_{j=1}^{\infty}2^{-j}x_{i,j})_{i\in\mathbb{N}}\;\;\;\;\;\;\;\;\big{(}x=(x_{i,j})_{i,j}\in\lambda^p(A)\big{)}
$$
is onto. Now, consider $Q^*:(\ell^p)^*=\ell^q\rightarrow\lambda^p(A)^*$, where $q=\frac{p}{p-1}$. By using \mbox{\cite[Corollary 27.24]{Meise}}, $rang(Q^*)$ is a closed subspace of $\lambda^p(A)^*$ which is non-reflexive. 
This example shows that not every closed subspace of a reflexive locally convex space is reflexive.

\begin{dif}
Let $\mathcal{A}$ be a Banach space and $B$ be a closed subspace of $\mathcal{A}$.
Following \mbox{\cite[page 293]{Meise}}, $\mathcal{A}$ is finite dimensional if and only if it is a Montel space. 
Therefore, $B$ and $\mathcal{A}/B$ are also Montel spaces. 
Clearly, for Fr\'echet spaces, every closed subspace of a Montel space is a Montel space.
But, Montel property are in general not inherited by quotients. For instance, consider the montel space $\lambda^1(A)$ in \mbox{D \ref{d5}} and the space $\ell^1$ which is not Montel. $\Box$
\end{dif}

Before continuing, we must recall that for a continuous linear map $T:\mathcal{A}\rightarrow\mathcal{B}$ between locally convex spaces $\mathcal{A}$ and $\mathcal{B}$, there is a unique map $\bar{T}$ in $\mathcal{L}(\mathcal{A}/ker(T),\mathcal{B})$ such that $T=\bar{T}\circ q$, where $q:\mathcal{A}\rightarrow\mathcal{A}/ker(T)$ is the quotient map \cite[Proposition 22.11]{Meise}.
Now, consider $rang(T)$ with the relative topology induced by $\mathcal{B}$. The map $T$ is said to be a topological homomorphism, if $\bar{T}$ is an isomorphism between $\mathcal{A}/ker(T)$ and $rang(T)$. 

\begin{dif}
Consider the short sequence of Fr\'echet spaces with continuous linear maps
\begin{equation}\label{e111}
\{0\}\rightarrow\mathcal{A}\overset{T}{\rightarrow}\mathcal{B}\overset{S}{\rightarrow}\mathcal{C}\rightarrow \{0\}.
\end{equation}
By applying \cite[Theorem 26.3]{Meise}, if the short sequence (\ref{e111}) is exact, then it is already topologically exact, that is $T$ and $S$ are topological homomorphisms. Also, by \cite[Proposition 26.4]{Meise}, the short sequence (\ref{e111}) is exact if and only if 
\begin{equation}\label{e222}
\{0\}\rightarrow\mathcal{C}^*\overset{S^*}{\rightarrow}\mathcal{B}^*\overset{T^*}{\rightarrow}\mathcal{A}^*\rightarrow \{0\}
\end{equation}
is exact. However, in general the short sequence (\ref{e222}) is not topologically exact.
Consider the K\"othe matrix $A$ and $Q\in\mathcal{L}(\lambda^p(A),\ell^p)$ defined in \mbox{D \ref{d5}}. By \mbox{\cite[Remark 27.23]{Meise}}, the dual sequence of the short exact sequence
$$
\{0\}\rightarrow ker(Q)\hookrightarrow\lambda^p(A)\overset{Q}{\rightarrow}\ell^p\rightarrow\{0\} 
$$
of Fr\'echet spaces is not topologically exact. \\
{\bf Remark.}
Let $\mathcal{A}$ be a Fr\'echet space and $B$ be a closed subspace of $\mathcal{A}$.
Consider the short exact sequence 
$$
\{0\}\rightarrow B \overset{j}{\rightarrow} \mathcal{A}\overset{q}{\rightarrow} \mathcal{A}/B \rightarrow \{0\},
$$
where $j:B\rightarrow\mathcal{A}$ is the inclusion and $q:\mathcal{A}\rightarrow\mathcal{A}/B$
is the quotient map. If the dual sequence 
$$
\{0\}\rightarrow(\mathcal{A}/B)^*\overset{q^*}{\rightarrow}\mathcal{A}^*\overset{j^*}{\rightarrow}B^*\rightarrow \{0\},
$$
is topologically exact, then we obtain the canonical isomorphisms
$B^*\cong\mathcal{A}^*/B^{\circ}$ and $(\mathcal{A}/B)^*\cong B^{\circ}$; see \cite[Remark 26.5]{Meise}.
$\Box$
\end{dif}

In \mbox{S \ref{spectrum}}, we mentioned an important result concerning the spectrum of Banach and Fr\'echet algebras. In the following, some other results in this field have been presented.

\begin{dif}
Let $\mathcal{B}$ be a commutative unital Banach algebra. Then, $\Delta(\mathcal{B})$
is a nonempty compact space. Following \cite[3.2]{Goldmann}, in general the spectrum of a commutative Fr\'echet algebra $\mathcal{A}$ is not a compact space. But, it is a nonempty hemicompact space. In fact, there exists a countable compact exhaustion $(K_n)_n$ of $\mathcal{A}$ such that for every compact subset $K\subseteq\mathcal{A}$ there is $n\in\mathbb{N}$ so that \mbox{$K\subseteq K_n$. $\Box$} 
\end{dif}

Before the next result, we recall that an ideal $I$ of a Fr\'echet algebra $\mathcal{A}$ is called maximal if $I\not=\mathcal{A}$ and $I$ is contained in no other proper ideal of $\mathcal{A}$.

\begin{dif}
Every maximal ideal in a commutative Banach algebra $\mathcal{B}$ is closed. Hence, there exists a one-to-one correspondence between the maximal ideals of $\mathcal{B}$ and the elements of $\Delta(\mathcal{B})$. Opposite to this result maximal ideals of commutative  Fr\'echet algebras are in general not closed. For example, $C(\mathbb{R})$
has some maximal ideal which is not closed; see \cite[Example 3.2.11]{Goldmann}. $\Box$
\end{dif}

\begin{dif}
Recall that a Hausdorff space is called a k-space if every subset intersecting each compact subset in a closed set is itself closed. Following \cite[chapter 3]{Goldmann}, examples of k-spaces are locally compact and first countable spaces. Also, the spectrum of a commutative Banach algebra is a k-space. But, there are Fr\'echet algebras whose spectrum is not a k-space; see \cite[Example 7.1.2]{Goldmann}\mbox{. $\Box$}
\end{dif}

Suppose that $\mathcal{A}$ is a commutative Fr\'echet algebra. Let us recall that a point $\varphi\in\Delta(\mathcal{A})$ is said to be a local peak point for $\mathcal{A}$ if there exist a neighborhood $U$ of $\varphi$ in $\Delta(\mathcal{A})$ and an element $a\in\mathcal{A}$ such that $\hat{a}(\varphi)=1$ and
$$
|\hat{a}(\psi)|<1\;\;\;\;\;\;\;\;(\psi\in U\setminus\{\varphi\}).
$$
Moreover, $a$ is called a peak point for $\mathcal{A}$ if $U$ can be taken to be $\Delta(\mathcal{A})$.

\begin{dif}
In the class of Banach algebras, by applying the Rossi's maximum principle, one can see that every locall peak point is a peak point; see \cite[Theorem 1.4.9]{Goldmann}. However, this result is not valid for Fr\'echet algebras; see also \mbox{\cite[Example 9.2.2]{Goldmann}}. Therefore, the Rossi's maximum modulus theorem fails to be true for Fr\'echet algebras. $\Box$
\end{dif}

\begin{dif}
Let $\mathcal{B}$ be a commutative unital Banach algebra and let $\mathcal{B}^{-1}$ denote the set of all invertible elements of $\mathcal{B}$. Due to the Arens-Royden theorem, for every $f\in C(\Delta(\mathcal{B}))^{-1}$ there exist $b\in\mathcal{B}^{-1}$ and $g\in C(\Delta(\mathcal{B}))$ such that $f=\hat{b}\;exp(g)$, where $exp(g)=\sum_{n=0}^{\infty}g^n/n!$; see \cite[Theorem 1.4.8]{Goldmann}.
It is an open question, whether the Arens-Royden theorem can be extended to Fr\'echet algebras. However, following \cite{Goldmann}, this theorem was proved for a special class of Fr\'echet algebras. $\Box$
\end{dif}

\section{\bf Some important generalizations in harmonic analysis}

In this section, we collect some important generalizations to the Fr\'echet algebras which have been obtained in recent years. We mark these results by the letter ``G", and we use the symbol ``$\Box$" after each one of them.

Let $\mathcal{B}$ be a Banach algebra. We recall that a left Banach $\mathcal{B}$-module $F$ is called flat if for every exact admissible sequence
$$
\{0\}\rightarrow E_1\rightarrow E_2\rightarrow E_3\rightarrow \{0\}
$$  
of right Banach $\mathcal{B}$-modules, the sequence
$$
\{0\}\rightarrow E_1\hat{\otimes}_{\mathcal{B}}F\rightarrow E_2\hat{\otimes}_{\mathcal{B}}F\rightarrow E_3\hat{\otimes}_{\mathcal{B}}F\rightarrow \{0\}
$$  
is exact. Also, $\mathcal{B}$ is amenable if for each Banach $\mathcal{B}$-bimodule $E$ every continuous derivation from $\mathcal{B}$ to the dual bimodule $E^*$ is inner. This definition was introduced by Johnson in \cite{Johnson1}. Later, Helemskii and  Sheinberg \cite{Hel0} showed that a Banach algebra $\mathcal{B}$ is amenable if and only if $\mathcal{B}^{\sharp}$, the unitization of $\mathcal{B}$, is a flat Banach $\mathcal{B}$-bimodule. In this case, $\mathcal{B}$ has a b.a.i. \mbox{\cite[Proposition 2.2.1]{Runde}}. In recent years, a large number of papers concerning amenability of Banach algebras have appeared; see for instance \cite{AZR,Alaghmandan,BR,Gordji2,Gordji3,Habibian,Medghalchi1}.
In \cite{Pir}, the arguments above have generalized as follows:

\begin{g}
Due to Pirkovskii \cite{Pir}, a Fr\'echet algebra $\mathcal{A}$ is amenable (resp. biflat) if $\mathcal{A}^{\sharp}$ (resp. $\mathcal{A}$) 
is a flat Fr\'echet $\mathcal{A}$-bimodule. Also, we say that $\mathcal{A}$ has a locally b.a.i. if for each zero neighborhood $U\subseteq\mathcal{A}$ there exists $C>0$ such that for each finite subset $F\subseteq\mathcal{A}$ and each $\varepsilon>0$ there exists $b\in CU$ with $a-ab\in\varepsilon U$ and $a-ba\in\varepsilon U$ for all $a\in F$. Note that if $\mathcal{A}$ is normable, then the notions of ``b.a.i.'' and ``locally b.a.i.'' are equivalent \mbox{\cite[Remark 6.4]{Pir}}. 
By applying \mbox{\cite[Lemma 9.4]{Pir}}, a Fr\'echet algebra $\mathcal{A}$ with a locally b.a.i. is amenable if and only if $\mathcal{A}$ is biflat. Furthermore, $\mathcal{A}$ is amenable if and only if it is isomorphic to a reduced projective limit of a sequence of amenable Banach algebras \mbox{\cite[Theorem 9.5]{Pir}}.
Let us note the following theorem proved by A. Pirkovskii \mbox{\cite[Theorem 9.6]{Pir}}.\\
{\bf Theorem.}
Let $\mathcal{A}$ be a Fr\'echet algebra. The following conditions are equivalent:
\begin{enumerate}
\item[(i)]
$\mathcal{A}$ is amenable;
\item[(ii)]
for each Banach $\mathcal{A}$-bimodule $E$, every continuous derivation from $\mathcal{A}$ to  $E^*$ is inner;
\item[(iii)]
for each Fr\'echet $\mathcal{A}$-bimodule $E$, every continuous derivation from $\mathcal{A}$ to  $E^*$ is inner.
\end{enumerate}
Pirkovskii also introduced the concepts of ``virtual diagonal", ``bounded approximate diagonal", and ``locally bounded approximate diagonal" in the class of Fr\'echet algebras. He then extended Johnson's theorem $($see \cite[Lemma 1.2]{Johnson} and \cite[Theorem 1.3]{Johnson}; see also \cite[Theorem 2.2.4]{Runde}$)$. Indeed, for a Fr\'echet algebra $\mathcal{A}$, the following statements are equivalent:
\begin{enumerate}
\item[(i)]
$\mathcal{A}$ is amenable;
\item[(ii)]
$\mathcal{A}$ has a locally bounded approximate diagonal.
\end{enumerate}
If, in addition, $\mathcal{A}$ is quasinormable, then (i) and (ii) are equivalent to the following:
\begin{enumerate}
\item[(iii)]
$\mathcal{A}$ has a bounded approximate diagonal;
\item[(iv)]
$\mathcal{A}$ has a virtual diagonal.
\end{enumerate}
See \cite[Theorem 9.7]{Pir}. $\Box$
\end{g}

\begin{g}
Ghahramani and Loy \cite{GhL} introduced and studied approximate amenability and contractibility of Banach algebras.
These notions were investigated in the class of Fr\'echet algebras by Lawson and Read \cite{LR}. 
A Fr\'echet algebra $\mathcal{A}$ is called approximately contractible if given any $\mathcal{A}$-bimodule $E$, and any continuous derivation $D:\mathcal{A}\rightarrow E$, there is a net $(x_{\alpha})\subseteq E$ such that
$$
D(a)=\lim_{\alpha}(a\cdot x_{\alpha}-x_{\alpha}\cdot a)\;\;\;\;\;\;\;\;\;\;\;(a\in\mathcal{A}).
$$
Every continuous derivation $D$ with this property is said to be approximately inner. Also, we recall that a Fr\'echet algebra $\mathcal{A}$ is approximately amenable if given any $\mathcal{A}$-bimodule $E$, every continuous derivation $D:\mathcal{A}\rightarrow E^*$ is approximately inner.
The authors in \cite{LR} generalized the important results of Banach algebras in this field. Following \cite{GhL}, all known approximately amenable Banach algebras have b.a.i. It is nonetheless hoped that it would yield Banach algebras without b.a.i. which had a form of amenability.
Note that there are examples of Fr\'echet algebras which are approximately contractible, but which do not have a b.a.i. For a good many Fr\'echet algebras without b.a.i., it was obtained either that the algebra is approximately amenable, or it is ``obviously" not approximately amenable because it has continuous point derivations; see \cite[section 3]{LR}.
Therefore, the situation for Fr\'echet algebras is quite close to what was hoped for Banach algebras. $\Box$
\end{g}

Weakly amenable Banach algebras were first introduced by Bade, Curtis and Dales in \cite{Bade}. The general definition is due to Johnson \cite{Johnson1}. Later, Dales, Ghahramani and \mbox{Gr$\o$nb$\ae$k} \cite{DGG} introduced and studied the notion of $n$-weak amenability of Banach algebras. Indeed, a Banach algebra $\mathcal{B}$ is $n$-weakly amenable if $H^1(\mathcal{B},\mathcal{B}^{(n)})=\{0\}$, where $\mathcal{B}^{(n)}$ is the $n$-th dual space of $\mathcal{B}$. In particular, $\mathcal{B}$ is weakly amenable if it is 1-weakly amenable. Recently, several papers dealing with this subject have appeared; see for example \cite{Gordji2,Ramezanpour,Samea}.

\begin{g}
Let $\mathcal{A}$ be a Fr\'echet algebra. The notion of weak amenability of $\mathcal{A}$ was studied by Rejali et al. \cite{ARR4}. In fact, the authors showed that some results in the field of weak amenability of Banach algebras can be generalized for the Fr\'echet algebras. For example, if $\mathcal{A}$ is weakly amenable, then $\mathcal{A}$ is essential \mbox{\cite[Theorem 2.3]{ARR4}.} Moreover, if $I$ is a quasinormable closed ideal in the Fr\'echet algebra $\mathcal{A}$ such that $\mathcal{A}/I$ and $I$ are weakly amenable, then $\mathcal{A}$ is weakly amenable, as well; see \mbox{\cite[Theorem 2.6]{ARR4}. $\Box$}
\end{g}

\begin{g}
Let $\mathcal{B}$ be a Banach algebra and $I$ be a closed ideal of $\mathcal{B}$.
The notions of $I$-weak amenability and ideal amenability of $\mathcal{B}$ was introduced and studied by Gorgi and Yazdanpanah in \cite{Gordji}. Indeed, $\mathcal{B}$ is $I$-weakly amenable if every continuous derivation $D:\mathcal{B}\rightarrow I^*$ is inner. Moreover, $\mathcal{B}$ is
called ideally amenable if $\mathcal{B}$ is $I$-weakly amenable for every closed ideal $I$ of $\mathcal{B}$. Clearly, every amenable Banach algebra is ideally amenable, and every ideally amenable Banach algebra is weakly amenable. However, the concepts of $I$-weak amenability and ideal amenability are different from amenability and weak amenability of Banach algebras (ibid.).

In \cite{RR}, Ranjbari and Rejali generalized these concepts and their results to the Fr\'echet algebras.
Later, the concept of vector-valued Lipschitz algebras in the class of Fr\'echet algebras was studied by the same authors in \cite{RR0}. Especially, they Proved that Lemma 2.1 of \cite{BR} holds for Fr\'echet algebras under a certain condition. Then, they studied ideal amenability of vector-valued Fr\'echet Lipschitz algebras; see \mbox{\cite[section 4]{RR0}} for details. $\Box$
\end{g}

Let $\mathcal{A}$ be a Banach algebra with a b.a.i. which is contained as closed ideal in a Banach algebra $\mathcal{B}$. Consider the pseudo-unital Banach $\mathcal{A}$-bimodule $E$, that is
$$
E=\mathcal{A}\cdot E\cdot\mathcal{A}=\{a\cdot x\cdot b:a,b\in\mathcal{A},\;x\in E\}.
$$
By using \cite[Proposition 2.1.6]{Runde}, $E$ is a Banach $\mathcal{B}$-bimodule in a canonical fachion. Moreover, for each derivation $D:\mathcal{A}\rightarrow E^*$, there is a unique extension derivation $\widetilde{D}:\mathcal{B}\rightarrow E^*$ of $D$ such that $\widetilde{D}$ is continuous with respect to the strict topology on $\mathcal{B}$ and the weak$^*$-topology on $E^*$. In \cite{AA}, we studied the strict topology on a Fr\'echet algebra and then we obtained a similar result as above.

\begin{g}
 Let $\mathcal{A}$ be a closed ideal in a Fr\'echet algebra $(\mathcal{B},p_{\ell})_{\ell\in\mathbb{N}}$, satisfying
$$
\big{\{}b\in\mathcal{B}: b\cdot\mathcal{A}=\mathcal{A}\cdot b=\{0\}\big{\}}=\{0\},
$$
and suppose that $\mathcal{B}$ has a b.a.i. with elements in $\mathcal{A}$.
Let $(E,q_n)_{n\in\mathbb{N}}$ be an arbitrary pseudo-unital Fr\'echet $\mathcal{A}$-bimodule.
Then, for every continuous derivation \mbox{$D:\mathcal{A}\rightarrow E^*$,} with respect to the topology induced by $(p_{\ell})_{\ell\in\mathbb{N}}$ on $\mathcal{A}$ and the strong topology on $E^*$, there exists a unique derivation $\widetilde{D}:\mathcal{B}\rightarrow E^*$ such that $\widetilde{D}$ is continuous with respect to the strict topology on $\mathcal{B}$ and the weak$^*$-topology on $E^*$; see \mbox{\cite[Theorem 3.2]{AA}.} $\Box$
\end{g}

\begin{g}
Following \cite{ARR2} and also \cite{ARR3}, a
Fr\'echet algebra $(\mathcal{A},p_{\ell})_{\ell\in\mathbb{N}}$ is called a Segal
Fr\'echet algebra in the Fr\'echet algebra $(\mathcal{B},q_{n})_{n\in\mathbb{N}}$
if the following conditions hold:
\begin{enumerate}
\item[(1)] 
$\mathcal{A}$ is a dense ideal in
$(\mathcal{B},q_{n})_n$;
\item[(2)] 
the identity map
$(\mathcal{A},p_{\ell})_{\ell}\rightarrow (\mathcal{B},q_{n})_n$ is
continuous; 
\item[(3)] 
the map $(\mathcal{B},q_{n})_n\times
(\mathcal{A},p_{\ell})_{\ell}\rightarrow (\mathcal{A},p_{\ell})_{\ell}$ with
$(a,b)\mapsto ab$, for all $a,b\in\mathcal{A}$, is jointly
continuous.
\end{enumerate}
This definition is in fact a generalization of the definition
abstract Segal algebras, in the Banach case; see \cite{Burn}. Also
in \cite[Definition 14]{B}, this definition has been introduced in
the field of topological algebras, as follows:

A topological algebra $(\mathcal A,\tau_{\mathcal A})$ is called a
$(\mathcal B,\tau_{\mathcal B})$-Segal algebra in the topological
algebra $(\mathcal B,\tau_{\mathcal B})$ if the following
assertions hold:
\begin{enumerate}
\item[(1)$^\prime$] 
$\mathcal A$ is a dense left ideal in
$(\mathcal B,\tau_{\mathcal B})$; 
\item[(2)$^\prime$] 
the identity
map from $(\mathcal A,\tau_{\mathcal A})$ into $(\mathcal
A,\tau_{\mathcal B})$ is continuous; 
\item[(3)$^\prime$] 
the multiplication
$$
(\mathcal B,\tau_{\mathcal B})\times(\mathcal A,\tau_{\mathcal
A})\rightarrow (\mathcal A,\tau_{\mathcal
A}),\;\;((a,b)\mapsto ab),
$$
is continuous.
\end{enumerate}

In \cite{ARR2}, Rejali et al. investigated Theorem 1.2 of \cite{Burn}
for Segal Fr\'echet algebras and obtained an analogous result. Indeed, if $(\mathcal{A},p_{\ell})_{\ell\in\mathbb{N}}$ is a proper Segal Fr\'echet algebra in $(\mathcal{B},q_{n})_{n\in\mathbb{N}}$ which contains an a.i. $(e_{i})_{i\in I}$, then $(e_{i})_{i\in I}$ cannot be bounded in $(\mathcal{A},p_{\ell})_{\ell\in\mathbb{N}}$ \cite[Proposition 3.5]{ARR2}. Furthermore, the authors characterized closed ideals of Segal Fr\'echet algebras, and showed that the ideal theorem is also valid for Fr\'echet algebras. In other words, every closed ideal of any Segal Fr\'echet algebra $(\mathcal{A},p_{\ell})_{\ell\in\mathbb{N}}$ in $(\mathcal{B},q_{n})_{n\in\mathbb{N}}$ is the intersection of a closed ideal of $\mathcal{B}$ with $\mathcal{A}$ \cite[Theorem 3.8]{ARR2}. $\Box$
\end{g}

Suppose that $\mathcal{A}$ is a commutative Fr\'echet algebra and $M(\mathcal{A})$ is the multiplier algebra with respect to the strict topology. As mentioned in section \ref{s0}, $\mathcal{A}$ is a $M(\mathcal{A})$-Segal algebra; see \cite[section 4]{AA} for details. 
In \cite{AA1}, we studied the notion of multiplier algebra of Fr\'echet algebras.

Let $\mathcal{B}$ be a Banach algebra and $\varphi\in\Delta(\mathcal{B})$. Kaniuth, Lau and Pym \cite{Kaniuth} introduced and studied the concept of (right) $\varphi$-amenability for Banach algebras as a generalization of left amenability of Lau algebras. 
We recall that, $\mathcal{B}$ is called (right) $\varphi$-amenable if there exists a continuous linear functional $m$ on $\mathcal{B}^*$ such that $m(\varphi)=1$ and $m(f\cdot a)=\varphi(a)m(f)$ for all $a\in\mathcal{B}$ and $f\in\mathcal{B}^*$. 
Furthermore, the notion of right character amenability of Banach algebras was introduced and studied by Monfared \cite{Monfared}. Indeed, a Banach algebra $\mathcal{B}$ is called right character amenable if for each $\varphi\in\Delta(\mathcal{B})\cup\{0\}$ and all Banach $\mathcal{B}$-bimodules $E$ with the left module action $a\cdot x=\varphi(a)x$, each continuous derivation $D:\mathcal{B}\rightarrow E^*$ is inner. 
Alaghmandan, Nasr-Isfahani, and Nemati \cite{Alaghmandan} characterized the character amenability of abstract Segal algebras. Specifically, they showed that $\varphi$-amenability of $\mathcal{B}$ is equivalent to 
$\varphi|_{\mathcal{A}}$-amenability of $\mathcal{A}$, where $\mathcal{A}$ is an abstract Segal algebra in $\mathcal{B}$. 
Right $\varphi$-amenability and right character amenability for Fr\'echet algebras were introduced and studied in \cite{ARR1}.
In the following, some related results about Segal Fr\'echet algebras are provided.

\begin{g}
Let $(\mathcal{A},p_{\ell})_{\ell\in\mathbb{N}}$ be a Segal Fr\'echet algebra in a Fr\'echet algebra $(\mathcal{B},q_n)_{n\in\mathbb{N}}$. Precisely, similar to Lemma 2.2 of \cite{Alaghmandan}, 
$$\Delta(\mathcal{A})=\{\varphi|_{\mathcal{A}}:\varphi\in\Delta(\mathcal{B})\};$$ 
see \cite[Lemma 6.7]{ARR1}.
Also, if $\varphi\in\Delta(\mathcal{B})$, then $\mathcal{B}$ is right $\varphi$-amenable if and only if $\mathcal{A}$ is right $\varphi|_{\mathcal{A}}$-amenable \cite[Theorem 6.8]{ARR1}.
The following is now immediate from results above. \\
{\bf Corollary.}
Let $(\mathcal{B},q_n)_{n\in\mathbb{N}}$ be a Fr\'echet algebra and let $(\mathcal{A},p_{\ell})_{\ell\in\mathbb{N}}$ be a Segal Fr\'echet algebra in $\mathcal{B}$. Suppose that $\mathcal{B}$ is right character amenable. Then, $\mathcal{A}$ is right $\varphi$-amenable for each $\varphi\in\Delta(\mathcal{A})$; see \cite[Corollary 6.9]{ARR1}. $\Box$
\end{g}

In addition to the above, we joint with Abtahi \cite{AA0} generalized the notion of Connes amenability of Banach algebras. We then proved some of the previous available results about this notion for Fr\'echet algebras. Rahnama and Rejali \cite{RahnamaA,RahnamaB} also introduced and studied the concepts of 
amenability modulo an ideal of Fr\'echet algebras and locally bounded approximate diagonal modulo an ideal of Fr\'echet algebras. 
Finally, the notion of BSE algebras in the class of Fr\'echet algebras was investigated by Rejali et al. in \cite{AmiriAbtahi1,AmiriAbtahi2,AmiriRejali}.

Before continuing, let us note the following arguments. Let $\mathcal{B}$ be a Banach algebra.
We know that on the bidual space $\mathcal{B}^{**}$ of $\mathcal{B}$, there are two multiplications, called the first and second Arens products, which make $\mathcal{B}^{**}$ into a Banach algebra. Moreover, $\mathcal{B}$ is called Arens regular if these two products coincide on $\mathcal{B}^{**}$. There are many papers concerning Arens regularity.
Also, this notion was investigated by Rejali et al. in several papers; see for example \cite{Abtahi,Baker,Khodsiani,Rejali00,Rejali0,Rejali}.
Note that if $\mathcal{B}^{**}$ is Arens regular, then so is $\mathcal{B}$. However, Pym \cite{Pym1} gives an example of an Arens regular Banach algebra $\mathcal{B}$ such that $\mathcal{B}^{**}$ is not Arens regular. Also, look at the following example.\\
{\bf Example.}
Consider the Banach algebras $\mathcal{A}_n=\ell^1(\mathbb{Z},w_{\frac{1}{n}})$ $(n\in\mathbb{N})$, where $w_{\alpha}(t)=(1+|t|)^{\alpha}$ $(t\in\mathbb{Z})$ for each $\alpha>0$. Let
\begin{eqnarray*}
\mathcal{A}&=&c_0\text{-}\prod_{n\in\mathbb{N}}\mathcal{A}_n=\big{\{}(a_n)\in\prod_{n\in\mathbb{N}}\mathcal{A}_n:\|a_n\|_n\rightarrow0\big{\}}\;\;\text{and}\\
\mathcal{B}&=&\ell^{\infty}\text{-}\prod_{n\in\mathbb{N}}\mathcal{A}_n=\big{\{}(a_n)\in\prod_{n\in\mathbb{N}}\mathcal{A}_n:\sup_{n\in\mathbb{N}}\|a_n\|_n<\infty\big{\}}.
\end{eqnarray*}
The following statements hold.
\begin{enumerate}
\item[(i)]
by \cite[page 35]{Dales}, $\mathcal{A}$ is Arens regular.
\item[(ii)]
by \cite[Example 9.2]{Dales}, $\mathcal{B}$ is not Arens regular.
\item[(iii)]
$\mathcal{B}$ is a closed sualgebra of $\mathcal{A}^{**}$ and consequently $\mathcal{A}^{**}$ is not Arens regular.
\end{enumerate}

Let $G$ be a locally compact (topological) group with a weight function $w$ on it.
Rejali and Vishki \cite{Rejali} showed that the following statements are equivalent:
\begin{enumerate}
\item[(i)]
$L^1(G,w)^{**}$ is Arens regular;
\item[(ii)]
$L^1(G,w)$ is Arens regular;
\item[(iii)]
$L^1(G)$ is Arens regular or $\Omega:G\times G\rightarrow (0,1]$, defined by
$$
\Omega(x,y)=\frac{w(xy)}{w(x)w(y)}\;\;\;\;\;(x,y\in G),
$$
is $0$-cluster, that is 
$\lim_m\lim_n\Omega(x_m,y_n)=0=\lim_n\lim_m\Omega(x_m,y_n)$,
for sequences $(x_m)$ and $(y_n)$ in $G$ with distinct elements.
\end{enumerate}
Moreover, Baker and Rejali proved that for a semigroup $S$ with a weight function $w$ on it, 
$\ell^1(S,w)$ is Arens regular, whenever $\ell^1(S)$ is  \cite{Baker}. 
Now, the following question is natural. \\
{\bf Question.}
Let $S$ be a semigroup and $w$ be a weight function on $S$. Are the following statements equivalent?
\begin{enumerate}
\item[(i)]
$\ell^1(S,w)^{**}$ is Arens regular.
\item[(ii)]
$\ell^1(S,w)$ is Arens regular.
\item[(iii)]
$\ell^1(S)$ is Arens regular or $\Omega$ is $0$-cluster.
\end{enumerate}
Later, Rejali et al. showed that if the set of all idempotents of an inverse semigroup $S$ is finite, then $\ell^1(S)$ is Arens regular if and only if $S$ is finite; see \mbox{\cite[Theorem 3.8]{Abtahi}.} Also in this case, if $S$ admits a bounded below weight for which $\Omega$ is $0$-cluster, then $S$ is countable; see \cite[Theorem 3.4]{Khodsiani}. By arguments above, another question arises.\\
{\bf Question.}
Let $S$ be an inverse semigroup with finitely many idempotents and $w$ be a weight function on $S$. Are the following statements equivalent?
\begin{enumerate}
\item[(i)]
$\ell^1(S,w)^{**}$ is Arens regular.
\item[(ii)]
$\ell^1(S,w)$ is Arens regular.
\item[(iii)]
$S$ is finite or $\Omega$ is $0$-cluster.
\end{enumerate}
{\bf Remark.}
In \cite[Example 3.9]{Abtahi}, the authors showed that there exists an inverse semigroup $S$ with infinitely many idempotents such that $\ell^1(S)$ is  Arens regular. 

Arens regularity of Fr\'echet algebras was defined by Zivari-Kazempour \cite{Zivari}. Indeed, a Fr\'echet algebra $\mathcal{A}$ is called Arens regular if the first and second Arens products coincide on $\mathcal{A}^{**}$. Later, we \cite{AA50} obtained some results in this field.

\begin{g}
In 1966, Gulick \cite{Gulick} studied Arens products on the bidual space of locally multiplicatively-convex topological (lmc for short) algebras with a hypo-continuous multiplication. We recall that the product 
$\mathcal{A}\times\mathcal{A}\rightarrow\mathcal{A}$, on a lmc algebra $\mathcal{A}$, is left (resp. right) hypo-continuous if for each zero neighborhood $U$ and for each bounded set $B$ of $\mathcal{A}$ there exists a zero neighborhood $V$ such that $VB\subseteq U$ (resp. $BV\subseteq U$). The multiplication in $\mathcal{A}$ is called hypo-continuous if it is both left and right hypo-continuous. 

Now, let $\mathcal{A}$ be a Fr\'echet algebra. Clearly, $\mathcal{A}$ belongs to the class of lmc algebras. Moreover, by applying \cite[III.5.1]{Sch}, the multiplication in $\mathcal{A}$ is hypo-continuous. Therefore, due to Gulick \cite{Gulick} and \cite[Proposition 25.9]{Meise}, $\mathcal{A}^{**}$ is always a Fr\'echet algebra under the strong topology and Arens products on $\mathcal{A}^{**}$.

Furthermore, Gulick generalized a theorem of Civin and Yood \cite{CY} to the lmc algebras, especially Fr\'echet algebras. He introduced a bicommutative algebra, that is an lmc algebra whose bidual is commutative. Then, he showed that subalgebras of bicommutative algebras are also bicommutative \mbox{\cite[Theorem 5.3]{Gulick}. $\Box$}
\end{g}

\begin{g}
Following \cite{CY}, a Banach algebra $\mathcal{B}$ has a weak right identity if there exists a net $(e_{\alpha})_{\alpha\in\Lambda}$ in $\mathcal{B}$ and an $M>0$ such that for every $\alpha\in\Lambda$, $\|e_{\alpha}\|<M$ and 
$$
\lim_{\alpha}f(ae_{\alpha}-a)=0\;\;\;\;\;\;\;\;(a\in\mathcal{B},\;f\in\mathcal{B}^*).
$$ 
Civin and Yood proved that $\mathcal{B}$ has a weak right identity if and only if $\mathcal{B}^{**}$ has a right unit with respect to the first Arens product \cite[Lemma 3.8]{CY}.
There is a generalization of this theorem in the class of lmc algebras especially Fr\'echet algebras.
Zivari-Kazempour proved that an lmc algebra $\mathcal{A}$ has a right b.a.i. (resp. left b.a.i) if and only if $\mathcal{A}^{**}$ has a right (resp. left) unit with respect to the first (resp. second) Arens product \cite[Theorem 3.2]{Zivari}. $\Box$
\end{g}

{\bf Acknowledgment.} The authors would like to thank the Banach
algebra center of Excellence for Mathematics, University of
Isfahan.

\vspace{9mm}

{\footnotesize \noindent
 Z. Alimohammadi\\
  Department of Mathematics,
   University of Isfahan,
    Isfahan, Iran\\
    z.alimohammadi@sci.ui.ac.ir; z.alimohamadi62@yahoo.com\\

\noindent
 A. Rejali\\
  Department of Mathematics,
   University of Isfahan,
    Isfahan, Iran\\
    rejali@sci.ui.ac.ir\\

\end{document}